\documentclass[a4paper,10pt]{article}
\usepackage{latexsym,amsmath,amsthm,amssymb,amsfonts}
\usepackage{natbib}
\usepackage{graphicx}

\author{\textbf{Miguel Brozos-V\'azquez} \\
{\small Dept. of Geometry and Topology} \\
{\small Universidade de Santiago de Compostela} \and
\textbf{Marco Antonio Campo-Cabana} \\
{\small Dept. of Applied Mathematics} \\
{\small Universidade de Santiago de Compostela} \and
\textbf{Jos\'e Carlos D\'iaz-Ramos} \\
{\small Dept. of Mathematics} \\
{\small University College Cork} \and
\textbf{Julio Gonz\'alez-D\'\i{}az}\footnote{Corresponding author:
Kellogg School of Management (CMS-EMS),
Northwestern University.
5100 Leverone Hall, 2001 Sheridan Road.
Evanston, IL 60208-2014.
Phone number: +1 847-467-1745.
Fax: +1 847-491-2530.
E-mail address: julio@northwestern.edu}\\
{\small Kellogg School of Management} \\
{\small Northwestern University} }

\title{Ranking Participants in Tournaments\\
 by means of Rating Functions}
\date{}

\newtheorem{property}{}

\newtheorem{theorem}{Theorem}
\newtheorem{lemma}[theorem]{Lemma}
\newtheorem{corollary}[theorem]{Corollary}
\newtheorem{prop}[theorem]{Proposition}

\theoremstyle{remark}

\newtheorem{example}{Example}

\providecommand{\norm}[1]{\lVert#1\rVert}
\providecommand{\abs}[1]{\lvert#1\rvert}

\newcommand{\TM}{A} 
\newcommand{\MM}{M} 
\newcommand{\GM}{B} 
\newcommand{\DM}{D^\MM} 
\newcommand{\TS}{\sigma^\rv} 
\newcommand{\AS}{\bar \sigma^\rv} 
\newcommand{\MS}{\mathcal M} 
\newcommand{\rv}{r} 
\newcommand{\pv}{p} 
\newcommand{\sv}{s} 

\begin{document}
\maketitle

\begin{abstract}
In this paper we bring a novel approach to the theory of
tournament rankings. We combine two different theories that are
widely used to establish rankings of populations after a given
tournament. First, we use the statistical approach of paired
comparison analysis to define the performance of a player in a
natural way. Then, we determine a ranking (and rating) of the
players in the given tournament. Finally, we show, among other properties, that the new ranking
method is the unique one satisfying a natural
consistency requirement.
\end{abstract}

\section{Introduction}\label{scIntro}

When there is a competition among the members of a population, the
fundamental problem is to rank these members according to their
strength.\footnote{We usually refer to the members of our
population as contestants or players, but they may also be other
objects such as scientific journals, political options or products to be tested.} In
certain cases this confrontation takes the form of a tournament in
which contestants play against themselves obtaining a certain
score in each match. The aim is to determine a final ranking after
all the matches have been played. Because of the wide range of
applications of ranking theory, the latter problem has already
been widely studied in management science, economic theory and
statistics. The books by \cite{Moon:1968} and \cite{Laslier:1997}
and the paper \cite{Iqbal:1986} discuss several applications of
the theory of tournament rankings.

Formally, a ranking of a population $N$ is a complete, reflexive
and transitive relation on $N$. To fix notation, we use the word
rating when we have a cardinal ranking, that is, not only do we
have an ordering of the contestants, but also a measure of the
intensities of the differences among them.

In this paper, we define a measure of the strength exhibited by
the players of a given tournament that we call \emph{performance}. More
specifically, we assume that there is a distribution function that
governs the random process associated with the competitive
environment in which the given tournament takes place. Within this
setting, the performance of a player is naturally defined as a
function of his results and the strength of his opponents.
Then, we define a new ranking method, the \emph{recursive
performance}, and show that it is the unique one that is consistent with the notion
of performance.

The most natural attempt to associate a ranking to a given
tournament is to use the \emph{scores ranking}. That is, rank the
players according to their total scores. \cite{Rubinstein:1980}
provides an axiomatic characterization of the scores ranking.
Nonetheless, it is often the case that several players have the
same score and, hence, in most scenarios this ranking method does
not provide an ordering of the players. On the other hand, most
ranking methods determine the ranking of a player according to the
results obtained in the tournament and the strength of the
opponents the player has played against, being this last feature
missing when using the scores ranking.\footnote{Indeed, this is
the reason why the scores ranking is mainly used in round-robin
tournaments, where each player faces the same opponents.} One
widely used ranking method that takes the previous considerations
into account is the {minimum violations ranking}, where one
violation consists of two players whose relative ranking differs
from the one induced by their result against each other. This
ranking method is discussed, for instance, in \cite{Goddard:1983}
and \cite{Iqbal:1986}. However, \cite{Stob:1985} is quite critical
with it and sticks up for the statistical approach of paired
comparison analysis initiated in \cite{Zermelo:1929} for chess
tournaments and rediscovered by \cite{Bradley:1952}.

In \cite{Bradley:1952}, each player $i$ is assumed to have a
strength parameter $\rv_i$. It is also assumed that there exists a
distribution $F$ such that $F(\rv_i,\rv_j)$ is the probability
that $i$ beats $j$. The objective is to fix a distribution $F$
that properly fits the available data and then, using statistical
tools, calculate the most likely values of the strength parameters
$\rv_i$. Once these values are calculated, they may be used to
define a rating for each player. As compared with the minimum
violations ranking, \cite{Stob:1985} emphasizes that this approach endogenizes
the relevance that must be given to the different victories
(losses) of the players, a feature that we also consider very
desirable. Two classic references within this framework are
\cite{Kendall:1940} and \cite{David:1988}.

Finally, a third approach comes from economic theory, where it is
often the case that an axiomatic approach is taken to determine a
ranking of the population $N$. First, it is assumed that there is
a matrix containing the relevant information about the paired
results of the different players; this matrix is usually referred
to as the tournament matrix. Then, a ranking method is defined as
a function that ascribes a ranking to each tournament matrix.
Next, the properties of the different methods are studied.
Finally, a ranking method, whose properties are suitable for a
given competitive environment, is chosen. To deepen into the
economic literature on rankings refer to \cite{Rubinstein:1980},
\cite{Liebowitz:1984}, \cite{Amir:2002}, \cite{Volij:2004}, and
\cite{Volij:2005}.

In this paper we consider a competitive environment with an
associated function $F$ that describes the underlying random
process. This function is determined from the data of the
historical confrontations of the players in a population. Then, we
assume that we have one more tournament and we want to rank the
players of that tournament according to their results. Next, we
bring together two widely used ideas. On the one hand, our ranking
method, the \emph{recursive performance}, is defined using a
recursive formula that resembles the Liebowitz-Palmer method
\citep{Liebowitz:1984} and other similar methods
\citep{Volij:2004} studied in economic theory. On the other hand,
our recursive formula uses the rating function $F$, so basic in
the statistical approach to paired comparison analysis. Moreover,
we show one property of our ranking method that is crucial for its
applicability. Namely, our ranking method is robust with respect
to the estimation of the function $F$, that is, small changes in
the function $F$ lead to small changes in the proposed rating.
Finally, concerning the computation of our ranking method, we show
that it reduces to solving a linear system.

Within the literature in which this paper is enclosed, there is a
family of tournaments that has received special attention: the
round-robin tournaments. In these tournaments each player plays
exactly once against any other player. The minimum violations
ranking is essentially thought for this class of tournaments.
\cite{Stob:1985} showed that, under quite general assumptions, the
ranking methods defined using the approach of \cite{Zermelo:1929}
lead to the same ordering as the scores ranking. We briefly discuss the behavior of our approach in round-robin tournaments at the end of Section~\ref{scRecursivePerformance}.

As an immediate application, the results obtained in this paper
can be applied to define new tie-breaking rules for disciplines
such as chess and Othello.\footnote{Remarkably, ties are always
present in tournaments in which pairings are drawn following the
Swiss pairing system, which is, along with the round-robin system,
the most widely used in these disciplines.} In our opinion, these
new tie-breaking rules improve the existent ones (see Sections
\ref{scRecursivePerformance} and \ref{scProperties} for details).

We briefly describe the contents of this paper. In
Section~\ref{scTournaments} we present the concept of tournament
and comment on the assumptions used throughout this paper. In
Section~\ref{scPerformance} we formally define and discuss the
notion of performance. Section~\ref{scRecursivePerformance} is the
core of our study; we introduce the so-called recursive
performance and give examples of tournaments in which it might be
immediately applied. In Section~\ref{scProperties} we discuss two
properties of the recursive performance ranking method: robustness
and consistency. Finally, in Section~\ref{scMaths} we prove the
results presented in Sections~\ref{scRecursivePerformance} and~\ref{scProperties}.

\section{Basic Definitions: Tournaments}\label{scTournaments}

We have a competitive environment in which confrontations between
the different players of a population take place along time. For
such an environment, there is a \emph{rating function} $F$ that
accurately describes the probabilities associated with the
different results of each match between any two given players.
This rating function is such that, given two players $i$ and $j$
with strength parameters $\rv_i$ and $\rv_j$, the probability that
$i$ beats $j$ is $F(\rv_i,\rv_j)$. Thus,
$F(\rv_i,\rv_j)=1-F(\rv_j,\rv_i)$.\footnote{For instance, chess
and Othello use rating systems based on functions that have
already been widely tested.} We refer to the strength parameters
$\rv_i$ as \emph{ratings}.

We work within the \emph{linear paired comparison model}
\citep{David:1988}. More specifically, we assume that there is a
strictly increasing continuous distribution function $F_l:\mathbb
R\to (0,1)$ such that $F(\rv_1,\rv_2)=F_l(\rv_1-\rv_2)$, that is,
the result of a game between any two players depends only on their
rating difference. The probability that $i$ beats $j$ goes to $1$
as $\rv_i-\rv_j$ grows and the probability that $i$ beats $j$ is
positive  regardless of the rating difference. Moreover, since
$F(\rv_i,\rv_j)=1-F(\rv_j,\rv_i)$, $F_l$ is symmetrically
distributed about zero. Also, note that the function $F_l^{-1}$ is
well defined.

The result of a confrontation between two players $i$ and $j$ may
be not only a win or a loss but any pair $(a_i,a_j)$ with $a_i,a_j\geq
0$, $a_i+a_j=1$. Within this scenario we might interpret $F(\rv_i,\rv_j)$ as the
expected score of player $i$ when facing player $j$.

We denote by $\MS_{k\times l}$ the vector space of real $k\times
l$ matrices. A \emph{tournament} is a pair $(N,\TM)$, where
$N=\{1,\ldots, n\}$ is the set of players and $\TM\in \MS_{n\times
n}$ is the \emph{tournament matrix}. The matrix $\TM$ is
non-negative and its main diagonal consists of zeros. The entry
$\TM_{ij}$ contains the score achieved by player $i$ in his
confrontations against player $j$. Note that we do not restrict
the non-zero entries in our matrix $\TM$ to be natural numbers as
in \cite{Volij:2005}. To each tournament $(N,\TM)$ we can
naturally associate a (symmetric) \emph{matches matrix} $\MM(\TM):=\TM+\TM^t$, \emph{i.e.}, $\MM(\TM)_{ij}$ is
the number of matches between $i$ and $j$. For notational
simplicity, when no confusion arises we denote $\MM(\TM)$ by
$\MM$. Since the $n$ players participate in the tournament, each
row of $\MM$ must have a nonzero entry. For each player $i$, let
$m_i:=\sum_{j=1}^n \MM_{ij}$ denote the total number of matches
played by $i$. Let $\rv\in \mathbb R^n$ be a vector of exogenously
given ratings of the players in the tournament $(N,\TM)$. We refer
to $\rv$ as the \emph{vector of initial ratings}.


Given a tournament $(N,\TM)$ and a vector of ratings $\rv$, we
define the \emph{total strength of the tournament} by
$\TS:=\sum_{i=1}^n m_i \rv_i $, that is, the contribution of each
player to the total strength of the tournament is weighted by the
number of matches he has played. Similarly, we define the
\emph{average strength of the tournament} by
$\AS:=\TS/\sum_{i=1}^n m_i$. Note that, if all the players have
played the same number of games, the average strength of the
tournament is just the average of the vector of initial ratings.

Summarizing, the primitives of our model are a tournament
$(N,\TM)$, a rating function $F$, and a vector of initial ratings
$\rv$. In this paper we present a ranking and rating method for
the tournament in question. Our ranking method reallocates among
the players the total strength of the tournament, $\TS$, with two
important features: first, the proposed ranking does not depend on
$\rv$ and, second, the difference between the proposed ratings for
any two players is also independent of $\rv$. That is, our method
is endogenous to $(N,\TM)$ and the vector $\rv$ is used just for
the sake of exposition.

Since the ranking method we define in this paper is anonymous, the
labels chosen for indexing the players are irrelevant. Based on
this fact, two tournaments that are equal up to labeling are said
equivalent. We make this definition precise. Denote by
$L_{\alpha\beta}\in \MS_{n\times n}$ the \emph{transposition
matrix} that swaps the $\alpha$th and $\beta$th entries of a
vector. A transposition matrix satisfies
$L_{\alpha\beta}^{-1}=L_{\alpha\beta}$ and, given 
$\GM\in \MS_{n\times n}$, the product $L_{\alpha\beta}\GM$ is the
same matrix $\GM$ but with rows $\alpha$ and $\beta$ interchanged.
Similarly, $\GM L_{\alpha\beta}$ interchanges columns $\alpha$ and
$\beta$ of $\GM$. The group $\Pi_n$ generated by the composition
of transposition matrices $L_{\alpha\beta}$ is isomorphic to the
group of permutations of $n$ elements. Given two populations $N$ and $N'$ with $n$ players, we say that two tournaments
$(N,\TM)$ and $(N',\TM')$ are \emph{equivalent} if there exists
$L\in\Pi_n$ such that $\TM=L\TM'L^t$. Note that the latter also
implies that $M=LM'L^t$. Since for each $L\in\Pi_n$ we have
$L^{-1}=L^t$, two equivalent tournaments have similar tournament
matrices and similar matches matrices.

A matrix $\GM\in\MS_{n\times n}$ is \emph{block diagonal},
respectively \emph{block anti-diagonal}, if
\[
\GM=\left(\begin{array}{c|c}*&0\\\hline 0&*\end{array}\right),\, \
\mbox{respectively}\ \GM=\left(\begin{array}{c|c}0&*\\\hline
*&0\end{array}\right).
\]

Note that, given a tournament $(N,\TM)$, $\TM$ is block diagonal
(anti-diagonal) if and only if $M$ is block diagonal
(anti-diagonal). We assume that our tournaments satisfy the
following two natural assumptions:

\begin{property}\label{A:AnonDiag}
The tournament $(N,\TM)$ is not equivalent to a tournament
$(N',\TM')$ such that $\TM'$ is block diagonal.
\end{property}

If the tournament matrix $\TM'$ is block diagonal, the tournament
has an internal division: there are two disjoint subsets of
players such that none of the players of one set has played
against anyone of the other set. This is a standard assumption in
the ranking's literature since each block may be considered as the
matrix of an independent tournament.

\begin{property}\label{A:AnonAntiDiag}
The tournament $(N,\TM)$ is not equivalent to a tournament
$(N',\TM')$ such that $\TM'$ is block anti-diagonal.
\end{property}

If $\TM'$ is block anti-diagonal, the tournament may be considered
as a team-tournament. There are two disjoint subsets (teams) such
that the players of each team have played only against the players
of the other, but not among themselves. Although similar to
\ref{A:AnonDiag}, this property is more subtle and has different
implications. In this case, in order to calculate the strength of
the players of one team, we need the strength of the players of
the other team that can only be calculated using the strength of
the players in the first team. This cyclic feature of
team-tournaments is the reason why, if \ref{A:AnonAntiDiag} is not
met, the iterative method we describe in
Section~\ref{scRecursivePerformance} does not necessarily
converge.

\section{The Notion of Performance}\label{scPerformance}

Let $(N,\TM)$ be a tournament. The \emph{vector of average scores}, $\sv$, is defined by
$\sv_i:=\sum_j A_{ij}/m_i$. Hereafter, we assume that $\sv \in (0,1)^n$. We define
$\DM:=\mathop{\rm diag}(m_1,\dots,m_n)\in \MS_{n\times n}$ and $\bar{\MM}:=(\DM)^{-1}\MM\in \MS_{n\times n}$, that is,
$\bar{\MM}_{ij}=\MM_{ij}/m_i$ is the number of confrontations
between $i$ and $j$ divided by the total number of matches played
by $i$.

Let $\rv\in \mathbb R^n$ be the vector of initial ratings and let
$F$ be the distribution function of the linear paired comparison model. The \emph{vector of $\rv$-performances}, $\pv^\rv \in \mathbb R^n$, is defined as
\[
\pv^{\rv}:= \bar \MM \rv+c, \quad\text{ where } c_i:=F_l^{-1}(\sv_i).
\]
Note that $(\bar \MM\rv)_i$ coincides with the average rating of
the opponents of player $i$. Hence, the $\rv$-performance of player $i$ is
the unique rating $p^\rv_i$ such that $F(\pv^\rv_i,(\bar \MM \rv)_i)=\sv_i$. Then, we might say that $\sv_i$ is $i$'s expected score against a player of rating $(\bar \MM \rv)_i$ if and only if $i$ has a rating $\pv^\rv_i$. This
justifies the name performance. 

The vector $\pv^\rv$ depends on the tournament and the rating
function associated with the competitive environment which the
tournament belongs to. Nevertheless, it also depends on the vector of initial ratings, which is exogenous to the tournament. On the
other hand, the $\rv$-performance of player $i$ takes into account
the strengths of his opponents; however we consider that also the
opponents of the opponents of $i$ should be used to calculate the
strengths of the opponents of $i$, and the $\rv$-performances do
not do so. We devote Section~\ref{scRecursivePerformance} to
address these two problems.

\begin{example}
The World Chess Federation (FIDE) has an official rating of
players called Elo. Elo's formula considers the distribution $F_l$
given by $F_l(\lambda)=1/\bigl(1+10^{-\frac{\lambda}{400}}\bigr)$. Hence, the
$\rv$-performance of player $i$ in a tournament is defined as $\pv^\rv_i=(\bar{\MM}\rv)_i-400\log_{10}(1/\sv_i-1)$, that is, the average of the Elos of his opponents plus a
correcting factor depending on the percentage of points achieved
by the player. Remarkably, this is one of the tie-breaking rules
recommended by the FIDE for chess tournaments.
\end{example}

The assumption $\sv \in (0,1)^n$ is needed in order to define the vector of performances correctly,
but this situation holds in most tournaments. Indeed, the $\rv$-performance ranking method is already used as a
tie breaking rule for chess tournaments, as we have just seen. The
idea of this method is to use the strengths of the opponents of
the players to define the rankings. Note that the vector $\rv$ is
the ``historical'' strength of the players whereas the vector $c$
is, essentially, the score of each player in the tournament.
Hence, $\pv^\rv$ measures the results of a player in relation
to the strength of his opponents.

\section{Recursive Performance}\label{scRecursivePerformance}

In the definition of $\rv$-performance, the vector $c$ depends crucially on
the rating function $F$ and, although our notation does not make
this dependence explicit, the rating function keeps being an
essential element of this paper.

The ranking associated with the vector of $\rv$-performances is
not a bad ranking for the tournament, but it heavily depends on
the initial ratings $\rv_1,\ldots,\rv_n$. The latter measure the
historical strength of the players, which might be different from
the strength exhibited by the players in the tournament. Moreover,
in the paired comparison literature, these ratings are often
calculated using the method of maximum likelihood, and thus, they
are subject to statistical errors.

The objective of this section and the next one is to
formally introduce a new ranking method, the \emph{recursive
performance}, and discuss some of its properties. Refer to Section
\ref{scMaths} for the technical results and their corresponding
proofs.

As we have already discussed in Section~\ref{scPerformance}, even
though $(\pv^r)_i$ is a better indicator of the strength of $i$ in
the given tournament than his initial rating, using the vector
$\pv^\rv$ as a ranking method has relevant shortcomings.
Nonetheless, it is natural to calculate a new performance by
replacing the initial ratings with the vector of performances. This
would have two main effects. Namely, it would shade the dependence
on $\rv$ and, given a player $i$, the new rating would take into
account, not only the opponents of $i$, but also their results
(\emph{i.e.}, using the opponents of the opponents of $i$ as
well). This suggests the iterative method $\pv^{(0)}:=\pv^\rv\,
(=\bar{\MM}\rv+c)$,
$\pv^{(l+1)}:=\pv^{\pv^{(l)}}\,(=\bar{\MM}\pv^{(l)}+c)$.
Unfortunately, this method does not necessarily converge. On the
other hand, the total strength of the tournament, $\sigma^\rv$, is not
preserved by the $\rv$-performance, that is, $\sum_i m_i \rv_i\neq \sum_i
m_i p^\rv_i$ in general. The inflation or deflation factor turns out
to be $\sum_i m_i c_i$. By distributing this factor among the
players, we define a new iterative method that does preserve the
total strength of the tournament.

Let $e\in\mathbb{R}^n$ be the vector $e=(1,\dots,1)$. Consider the
following rescaling of $c$,
\[
\hat{c}:=c-\left(\frac{\sum_i m_i c_i}{\sum_i m_i}\right)e,
\]
which we discuss below. We define the \emph{iterated performance}
as the iterative method
\[
\begin{array}{rcl}
\hat \pv^{(0)}&:=&\bar \MM \rv+\hat{c}\\[1ex]
\hat \pv^{(l)}&:=&\bar \MM \hat \pv^{(l-1)}+\hat{c},\quad
l\in\mathbb{N}.
\end{array}
\]
At each step $l$, this iterative method gives the same rating as the
previous one up to a constant proportional to $e$ and, hence,
\emph{the two proposed rankings are always the same}. This is
proved in the following lemma.

\begin{lemma}\label{order-not-changed}
For each $l\in\mathbb{N}$, $\pv^{(l)}-\hat
\pv^{(l)}=(l+1)\left(\frac{\sum_i m_i c_i}{\sum_i m_i}\right)e$.
\end{lemma}
\begin{proof}
Since by definition of $\bar{\MM}$, $\bar{\MM}e=e$, the result follows
by an induction argument.
\end{proof}

In Section \ref{scMaths} we prove that, under our assumptions, the
iterated performance converges. More specifically, we have the
following result.

\begin{theorem}\label{thConvergence}
Let $(N,\TM)$ be a tournament satisfying assumptions
\ref{A:AnonDiag} and \ref{A:AnonAntiDiag} and such that
$\bar\TM_{ij}\in (0,1)$ whenever $\MM_{ij}\neq 0$. Let
$\rv\in\mathbb{R}^n$ be a vector of initial ratings and $F_l$ the
distribution of the linear paired comparison model. Then the
iterated performance converges.
\end{theorem}

Motivated by Theorem~\ref{thConvergence}, we define the main concept of this paper, the \emph{recursive performance}, as the limit of the iterated performance:
$\hat{\pv}:=\lim_{l\to\infty}\hat{\pv}^{(l)}$. Taking limits in the equality
$\hat{\pv}^{(l)}=\bar \MM \hat \pv^{(l-1)}+\hat{c}$, we get that
$\hat{\pv}$ is a solution of the linear system
\begin{equation}\label{linear-system:performance}
(I-\bar{\MM})x=\hat c,
\end{equation}
where $I\in \MS_{n\times n}$ is the identity matrix. If
\ref{A:AnonDiag} holds, by
Theorem~\ref{teor:bound-eigenvalues:generalized} (iii) below the
matrix $I-\bar{\MM}$ has rank $n-1$. Then, since $\bar{\MM}e=e$, the
whole set of solutions of \eqref{linear-system:performance} is
given by $\hat{\pv}+\text{span}\{e\}$. The different solutions of
\eqref{linear-system:performance} arise from different initial
vectors of ratings $\rv$. It is important to note that \emph{all the
solutions propose the same ranking}.

Following the previous discussion, even if \ref{A:AnonAntiDiag}
does not hold, we can unambiguously associate a ranking to each
linear system \eqref{linear-system:performance} as far as
\ref{A:AnonDiag} is met.

\begin{example}
An \emph{ideal chess-like tournament} is a tournament in which all
the players play the same number of rounds, say $k$.\footnote{Most
tournaments in disciplines such as chess and Othello have this
property.} Thus, for an ideal chess-like
tournament $(N,\TM)$ we have $\DM=kI$ and, hence,
$\bar{\MM}=\MM/k$.  By Corollary \ref{thConservacion} below, $\sum_{i}m_i
\hat{\pv}_i^{(l)}=\sum_{i} m_i \rv_i=\TS$  for all $l\in \mathbb N$.
Since $m_i=k$ for all $i\in N$, in an ideal chess-like tournament we have $\sum_{i}
\hat \pv^{(l)}_i=\sum_{i} r_i$. This shows that, using $\hat{c}$
instead of $c$, we adjust the vectors after each iteration to
ensure that the sum of the ratings after each step remains
constant. The average strength of an ideal tournament, $\AS$,
coincides with the average of the components of $\rv$. In each
iteration the method proposes a way to divide the total strength
of the tournament, $\TS$, among the players. That is, by working
with $\hat{c}$ instead of $c$, we avoid inflation or deflation in
the vectors of iterated performances. Since the recursive
performance is the limit of such vectors, it also provides a way
of dividing the total strength of the tournament among the
players. In a general tournament, the same property holds, but in
this case the average strength of the tournament is calculated as
a weighted average (by the $m_i$'s) of the components of $\rv$.
\end{example}

\begin{example}
A \emph{round-robin tournament} $(N,\TM)$ is a tournament in which
all the players have played exactly once against each other. That is, $\MM_{ij}=1$ if $i\neq j$ and
${\MM}_{ii}=0$.\footnote{Round-robin tournaments have a special
structure and different approaches have been taken to define
ranking methods within this family of tournaments. See for
instance \cite{Daniels:1969}, \cite{Stob:1985} and Chapter~6.1 in
\cite{David:1988}.} Within this family of tournaments, the ranking proposed by the recursive
performance has an appealing feature: it coincides with the ranking
proposed by the vector of average scores (the standard scores ranking method). To see this, let $x$ be a solution of
the system \eqref{linear-system:performance}. Then, the claim follows from
the equality $\hat{c}_i-\hat{c}_j=x_i-\sum_{k\neq
i}\frac{x_k}{(n-1)}-x_j+\sum_{k\neq
j}\frac{x_k}{(n-1)}=\frac{n}{n-1}(x_i-x_j)$. This is not
surprising since the ranking proposed by the recursive performance
uses both the scores of the players and the scores of the
opponents, but all the players have the same
opponents. 
\end{example}

\section{Main Properties:  Robustness and Consistency}
\label{scProperties}

Since all the solutions of (\ref{linear-system:performance}) propose the same ranking, the recursive performance ranking method does not depend on the vector of initial ratings.
Moreover, the recursive performance ranking can be unambiguously
calculated for tournaments in which there are unrated players
(players with no historical results). Thus, if there is an unrated
player, we can assign him an arbitrarily chosen rating and this
election does not affect the final ranking. The recursive
performance proposes a way to divide the total strength of the
tournament among the players. Therefore, when used as a rating
method, there is neither inflation nor deflation with respect to
the initial ratings.


\subsection*{Robustness in $F$}

In order to apply the recursive performance in a given competitive
environment, we first need to estimate the function $F$ that
governs it. This estimation is subject to error and hence we need
to ensure that the ranking method is \emph{robust} in $F$, that
is, small changes of $F$ lead to small changes of the recursive
performance. This follows immediately from the fact that the
recursive performance is a solution of the linear system
(\ref{linear-system:performance}).

\subsection*{Consistency with $F$}

Given two vectors of ratings $\rv^1,\rv^2 \in \mathbb R^n$, we say
that they are \emph{essentially identical for function $F$} if
$F(\rv^1_i,\rv^1_j)=F(\rv^2_i,\rv^2_j)$ for all $i\neq j$. That
is, two vectors are essentially identical if they lead to the same
predictions for any given match between players in $N$. In linear
paired comparisons analysis this is equivalent to the existence of
$\lambda\in\mathbb{R}$ such that $r^1=r^2+\lambda e$.

\begin{lemma}
Let $(N,\TM)$ be a tournament and $F_l$ a linear rating function.
Then, all the solutions of (\ref{linear-system:performance}) are
essentially identical for function $F_l$.
\end{lemma}

\begin{proof}
It follows from the fact that all the solutions of
\eqref{linear-system:performance} are of the form
$\hat{p}+\text{span}\{e\}$.
\end{proof}

Because of this property, we make a slight abuse of language and
use the words recursive performance to encompass all the solutions
of (\ref{linear-system:performance}). Then, the total strength of
the tournament suffices to pin a unique rating.

Now, we are ready to introduce a consistency property. We say
that a rating $\rv\in \mathbb R^n$ is \emph{consistent with $F$
for tournament $(N,\TM)$} if the $\rv$-performance rating vector
is essentially identical to $\rv$. This consistency property can
also be extended to any ranking (not necessarily cardinal). A
ranking $\succeq$ is \emph{consistent with $F$} for the tournament
$(N,\TM)$ if there is a rating $\rv\in \mathbb R^n$ that is
consistent with $F$ and whose induced ranking coincides with
$\succeq$.

At the end of Section \ref{scMaths} we give a proof of the
following result.

\begin{prop}\label{thConsistency}
Let $(N,\TM)$ be a tournament and $F_l$ a linear rating function.
	Then, the ranking induced by the recursive performance is the
unique one that is consistent with $F_l$.
\end{prop}

\section{Mathematical results}\label{scMaths}

In this section we prove the technical results we have used
throughout Sections \ref{scRecursivePerformance} and
\ref{scProperties}. Although these results are stated for
tournaments, they may be written just in terms of linear algebra.
We follow \cite{Ciarlet:1989}.

A linear iterative method is (globally) convergent if and only if the
eigenvalues of the corresponding matrix are, in absolute value,
less than $1$. For any tournament $(N,\TM)$ we have $\bar{\MM}e=e$
and thus $1$ is an eigenvalue of $\bar{\MM}$. In this section we
prove that, under the assumptions \ref{A:AnonDiag} and
\ref{A:AnonAntiDiag}, the iterated performance restricts to a
vector subspace where the absolute values of the eigenvalues of
$\bar{\MM}$ are less than $1$ and, hence, the method converges.

Let $(N,\TM)$ be a tournament. For each $v,w\in\mathbb{R}^n$, we
define $\langle v,w\rangle:=v^t \DM w=\sum_{i=1}^n m_i v_i w_i$. Since $\DM$ is a diagonal
matrix and $m_i>0$ for all $i$, $\langle\cdot,\cdot\rangle$ is an
inner product in $\mathbb{R}^n$, which generalizes the Euclidean
inner product. The former, which depends on the tournament, is
referred to as the \emph{inner product associated with $(N,\TM)$}.

If $(N,\TM)$ is a tournament, then $\MM$ is a symmetric matrix but
$\bar {\MM}$ is not symmetric in general. However, there is a kind
of symmetry in $\bar {\MM}$, namely, $\bar {\MM}_{ij}=0$ if and
only if $\bar {\MM}_{ji}=0$. Motivated by this fact, we say that
two matrices $B\in \MS_{k\times l}$ and $C\in \MS_{l\times k}$ are
\emph{null-transpose} if for each $i \in \{1,\ldots,k\}$, and each
$j \in \{1,\ldots,l\}$, $B_{ij}=0$ if and only if $C_{ji}=0$. With
a slight abuse of notation we denote by $B^{nt}$ a matrix that is
null-transpose of $B$. Note that, although $B^{nt}$ is not unique,
$B^{nt}=0$ if and only if $B=0$.

\begin{theorem}\label{teor:bound-eigenvalues:generalized}
Let $(N,\TM)$ be a tournament and $\langle\cdot,\cdot\rangle$ its
associated  inner product. Then
\begin{enumerate}

\item The matrix $\bar {\MM}$ is self-adjoint with respect to
$\langle\cdot,\cdot\rangle$. Moreover, it is diagonalizable, its
eigenvalues are real and the eigenspaces are orthogonal with
respect to $\langle\cdot,\cdot\rangle$.

\item If $\lambda$ is an eigenvalue of $\bar {\MM}$, then
$\abs{\lambda}\leq 1$.

\item $(N,\TM)$ satisfies \ref{A:AnonDiag} if and only if the
multiplicity of the eigenvalue $1$ is one.

\item $(N,\TM)$ satisfies \ref{A:AnonAntiDiag} if and only if $-1$
is not an eigenvalue of $\bar {\MM}$.
\end{enumerate}
\end{theorem}

\begin{proof}
Let $(N,\TM)$ and $(N',\TM')$ be equivalent tournaments with
matches matrices $\MM$ and $\MM'$ respectively. Then, there exists
$L\in\Pi_n$ such that $\MM=L\MM'L^{-1}$. Clearly,
$\DM=LD^{M'}L^{-1}$ and thus $\bar{\MM}=L\bar{\MM}'L^{-1}$. Hence,
$\bar{\MM}$ and $\bar{\MM}'$ are similar matrices and their eigenvalue
structure is the same. Therefore, we may make, without loss of
generality, any assumption regarding the ordering of the indices.
We also recall that $\bar{\MM}_{ij}\geq 0$ for all
$i,j\in\{1,\dots,n\}$. For each $k \in \mathbb{N}$, we define $e^k:=(1,\dots,1)\in \mathbb{R}^k$. Note that
$e^n=e$.

Claim (i): Since $\MM$ and $\DM$ are symmetric, $\langle v,\bar
{\MM}w\rangle= v^t \DM\bar {\MM} w=v^t \MM w=w^t \MM v=w^t \DM\bar
{\MM}v=(\bar {\MM}v)^t\DM w=\langle \bar {\MM}v,w\rangle$. The
second part follows from the spectral theorem.

Claim (ii): The matrix norm $\norm\cdot_\infty$ is defined as
$\norm{B}_\infty:=\max_{i\in \{1,\ldots,n\}} \sum_{j=1}^n \lvert
B_{ij} \rvert$ for any $B\in\MS_{n\times n}$. By
definition we have $\norm{\bar {\MM}}_\infty=1$ and hence (ii)
follows from Theorem 1.4-3 in \cite{Ciarlet:1989}.

Claim (iii): Assume that $(N,\TM)$ does not
satisfy~\ref{A:AnonDiag}. Then, $\bar {\MM}$ may be written as
\[
\bar {\MM}=\left(\begin{array}{c|c}C&0\\ \hline
0&E\end{array}\right),\quad\mbox{with}\ C\in\MS_{k\times
k}\ \mbox{and}\ E\in\MS_{l\times l}.
\]
Hence, $\bar {\MM}(e^k|0)=(e^k|0)$ and $\bar
{\MM}(0|e^l)=(0|e^l)$. Thus, $1$ has multiplicity at least~$2$.

Conversely, assume that $1$ has multiplicity greater than one.
Then, since $\bar {\MM}$ is diagonalizable there exists $v\in \mathbb
R^n$, linearly independent from $e$, such that $\bar {\MM}v=v$.
Assume that $v_1=1$ and that the components of $v$ are
decreasingly ordered, that is, $v_i\geq v_j$ for $i>j$. Let $k\in
\mathbb N$ be such that $v_k=1$ and $v_{k+1}<1$. Since $v$ and $e$
are linearly independent, $k<n$ and $\bar {\MM}$ may be decomposed as
\begin{equation}\label{eqDecomposition}
\bar {\MM}=\left(\begin{array}{c|c}
C_1&E\\
\hline E^{nt}&C_2
\end{array}\right),\ \mbox{where}\ C_1\in\MS_{k\times k},\
C_2\in\MS_{(n-k)\times (n-k)}\ \mbox{and}\ E
\in\MS_{k\times (n-k)}.
\end{equation}
Now, if $E$ has a nonzero row, namely row $i$, we get
\[
1=v_i=(\bar {\MM}v)_i=\sum_{j=1}^k\bar {\MM}_{ij}
+\sum_{j=k+1}^n\bar {\MM}_{ij}v_j<\sum_{j=1}^n\bar {\MM}_{ij}=1,
\]
contradiction. This proves $E=0$ and $E^{nt}=0$, which is a
contradiction with \ref{A:AnonDiag}.

Claim (iv): Assume $(N,\TM)$ does not satisfy
\ref{A:AnonAntiDiag}. Then we may write $\bar {\MM}$ as
\[
\bar {\MM}=\left(\begin{array}{c|c}0&C\\
\hline E&0\end{array}\right), \quad\mbox{with}\
C\in\MS_{k\times l}\ \mbox{and}\ E\in\MS_{l\times
k}.
\]
Hence, $\bar {\MM}(e^k|-e^l)=-(e^k|-e^l)$ and $-1$ is an eigenvalue of
$\bar {\MM}$.

Conversely, assume that $-1$ is an eigenvalue of $\bar {\MM}$. Let
$v\in \mathbb R^n$ be such that $\bar {\MM}v=-v$ and $1=v_1\geq
v_2\geq\dots\geq v_n\geq -1$. Again, there exists $k\in \mathbb
N$, $k<n$ such that $v_k=1$ and $v_{k+1}<1$. Hence, $\bar {\MM}$ may
be decomposed as in \eqref{eqDecomposition}.

We show that $C_1=0$. Let $i\in\{1,\dots,k\}$. Since
${\sum_{j=1}^n\bar {\MM}_{ij}=1}$, we have
\[
-1=-v_i=(\bar {\MM}v)_i=\sum_{j=1}^k\bar {\MM}_{ij}
+\sum_{j=k+1}^n\bar {\MM}_{ij}v_j\geq\sum_{j=1}^k\bar {\MM}_{ij}
+\Bigl(1-\sum_{j=1}^k\bar {\MM}_{ij}\Bigr)(-1)=
2\sum_{j=1}^k\bar {\MM}_{ij}-1.
\]
Hence, $\sum_{j=1}^k\bar {\MM}_{ij}\leq 0$. Since $\bar {\MM}_{ij}\geq
0$, the $i$th row of $C_1$ is zero. Therefore, $C_1=0$.

Note that $v_n=-1$ since, otherwise, taking $i\in\{1,\dots,k\}$ we
get $-1=-v_i=(\bar {\MM}v)_i=\sum_{j=k+1}^n\bar {\MM}_{ij}v_j>-1$. Let
$l\in\mathbb{N}$ be such that $v_{n-l}>-1$ and $v_j=-1$ for $j\geq
n-l$. Clearly, $l>0$. Then, we may further decompose $\bar {\MM}$ as
\[
\bar {\MM}=\left(\begin{array}{c|c|c}
0&E_1&E_2\\
\hline E_1^{nt}&C_{21}&C_{22}\\
\hline E_2^{nt}&C_{22}^{nt}&C_{23}
\end{array}\right)\
\mbox{with}\ E=\left(E_1 \vert E_2\right),\
C_2=\left(\begin{array}{c|c}C_{21}&C_{22}\\\hline
C_{22}^{nt}&C_{23}\end{array}\right)\ \mbox{and}\
E_2\in\MS_{k\times l}.
\]
If $k+l=n$, then this second decomposition is trivial ($E=E_2$ and
$C_2=C_{23}$). Otherwise, we claim that $E_1=0$. If the $i$th row
of $E_1$ is nonzero, then there is $j\in\{k+1,\dots,n-l\}$ such
that $v_j>-1$. Hence,
$-1=-v_i=(\bar {\MM}v)_i=\sum_{j=k+1}^{n}\bar {\MM}_{ij}v_j>-1$,
contradiction. Therefore, $E_1=0$.

Now, we show that $(C_{22}^{nt}\vert C_{23})$ is zero. If the
$i$th row of $(C_{22}^{nt}\vert C_{23})$ is nonzero, then there is
$j\in\{k+1,\dots,n\}$ such that $v_j<1$. Hence,
$1=-v_i=(\bar {\MM}v)_i=\sum_{j=1}^n\bar {\MM}_{ij}v_j<1$,
contradiction. Therefore, the above decomposition reduces to
\[
\bar {\MM}=\left(\begin{array}{c|c|c}
0&0&E_2\\
\hline 0&C_{21}&0\\
\hline E_2^{nt}&0&0
\end{array}\right),
\]
which contradicts \ref{A:AnonAntiDiag}.\footnote{Note that, if the
latter decomposition happens to be nontrivial (presence of
$C_{21}$), we also violate \ref{A:AnonDiag}.}
\end{proof}

\begin{corollary}\label{corollary-Amatrix-convergence:generalized}
Let $(N,\TM)$ be a tournament satisfying assumptions
\ref{A:AnonDiag} and \ref{A:AnonAntiDiag}. Then,
\[
\lim_{l\to\infty}\bar {\MM}^l=\frac{1}{\sum_{i=1}^n m_i} \left(
{\begin{array}{ccc}
m_1&\cdots&m_n\\
\vdots&&\vdots\\
m_1&\cdots&m_n
\end{array}}\right).
\]
\end{corollary}

\begin{proof}
Let $v\in\mathbb{R}^n$. By Theorem
\ref{teor:bound-eigenvalues:generalized} there exists a basis of
eigenvectors $\{w_1=e,w_2,\dots,w_n\}$, orthogonal with respect to
the inner product associated with $(N,\TM)$, which we denote by
$\langle\cdot,\cdot\rangle$. For each $i\in\{1,\dots,n\}$, let
$\lambda_i$ be the eigenvalue associated with $w_i$. Now, it
suffices to show that $(\lim_{l\to\infty}\bar
{\MM}^l)v=(\frac{1}{\mathop{\rm tr}\DM}e e^t \DM)v$ for all $v\in
\mathbb R^n$. For each $i\in\{1,\dots,n\}$, we have
$\langle(\lim_{l\to\infty}\bar {\MM}^l)
v,w_i\rangle=\lim_{l\to\infty}\langle v,\bar {\MM}^l
w_i\rangle=(\lim_{l\to\infty}\lambda_i^l)\langle v,w_i\rangle$.
Since $\lambda_1=1$ and $\lvert\lambda_i\rvert<1$ for each
$i\in\{2,\dots,n\}$, we get $(\lim_{l\to\infty}\bar {\MM}^l)
v=(\langle v,e\rangle/\langle e,e\rangle)e$. By definition of
$\langle\cdot,\cdot\rangle$, we have $\langle
e,e\rangle=\mathop{\rm tr}\DM$ and $\langle v,e\rangle =e^t\DM v$,
from where the result follows.
\end{proof}

\begin{corollary}\label{thConservacion}
Let $(N,\TM)$ be a tournament satisfying \ref{A:AnonDiag} and
\ref{A:AnonAntiDiag}. Let $\langle\cdot,\cdot\rangle$ be its
associated~inner product. Then, for all $l\in\mathbb{N}$,
$\langle\hat{\pv}^{(l)},e\rangle=\langle \rv,e\rangle$, or
equivalently, $\sum_{i}m_i \hat{\pv}_i^{(l)}=\sum_{i} m_i
\rv_i=\TS$.
\end{corollary}

\begin{proof}
Recall that $\bar {\MM}$ is self-adjoint with respect to
$\langle\cdot,\cdot\rangle$ by Theorem
\ref{teor:bound-eigenvalues:generalized}~(i). By Corollary~\ref{corollary-Amatrix-convergence:generalized},
$\lim_{l\to\infty}\bar {\MM}^l$ exists and
$(\lim_{l\to\infty}\bar {\MM}^l)c=\bigl((\sum_i m_i c_i)/(\sum_i
m_i)\bigr)e$. Hence, $\langle\hat{c},e\rangle=\langle
c,e\rangle-\langle (\lim_{l\to\infty}\bar {\MM}^l)c, e\rangle=\langle
c,e\rangle-\lim_{l\to\infty}\langle c,\bar {\MM}^l e\rangle=0$. Then,
$\langle\hat{\pv}^{(0)},e\rangle =\langle
\rv,\bar {\MM}e\rangle+\langle\hat{c},e\rangle=\langle \rv,e\rangle$. By
induction, $\langle\hat{\pv}^{(l)},e\rangle =\langle
\hat{\pv}^{(l-1)},\bar {\MM}e\rangle+\langle\hat{c},e\rangle=\langle
\rv,e\rangle$.
\end{proof}

We are now ready to prove the main result of this paper.

\begin{proof}[Proof of Theorem \ref{thConvergence}]
Defining $q^{(l)}:=\hat \pv^{(l)}-\rv$ for $l\geq 0$ we have the
equivalent iterative method $q^{(0)}=\bar {\MM} \rv+\hat c-\rv$
and $q^{(l)}=\bar {\MM} q^{(l-1)}+q^{(0)}$, $l\in\mathbb{N}$. Let
$\langle\cdot,\cdot\rangle$ be the inner product associated with
$(N,\TM)$. By Corollary \ref{thConservacion}, $\langle
q^{(l)},e\rangle=0$ for all $l\in\mathbb{N}$. Therefore, the
iterative method restricts to the vector subspace $e^\perp$. By
Theorem~\ref{teor:bound-eigenvalues:generalized}, the absolute
values of the eigenvalues of $\bar {\MM}_{|e^\perp}$ are smaller
than $1$. Hence, the iterative method converges
\citep{Ciarlet:1989}.
\end{proof}

\begin{proof}[Proof of Proposition \ref{thConsistency}]
Let $\succeq$ be a ranking consistent with $F_l$. Then, there
exists a rating $r\in\mathbb{R}^n$ and a constant
$\lambda\in\mathbb{R}$ such that $r=p^r+\lambda e=\bar\MM
r+c+\lambda e$. Taking inner product with $e$ (as in Corollary
\ref{corollary-Amatrix-convergence:generalized}) we get $\langle
r,e\rangle=\langle r,e\rangle+\langle c,e\rangle+\lambda\langle
e,e\rangle$. Hence, $\lambda=-\langle c,e\rangle/\langle
e,e\rangle$ and thus $c+\lambda e=\hat{c}$. This implies that $r$
is a solution of \eqref{linear-system:performance} and the result
follows.
\end{proof}

\end{document}